\documentclass [11pt,oneside,a4paper,mathscr]{amsart}

\usepackage{amscd,amsmath,amssymb,euscript}
\usepackage[frame,cmtip,curve,arrow,matrix,line,graph]{xy}

\title{T-structures on some local Calabi-Yau varieties}
\author{Tom Bridgeland}

\date{}

\newtheorem{thm}{Theorem}[section]

\newtheorem{prop}[thm]{Proposition}
\newtheorem{lemma}[thm]{Lemma}
\newenvironment{pf}{\paragraph{Proof}}{\qed\par\medskip}
\theoremstyle{definition}
\newtheorem{defn}[thm]{Definition}

\newtheorem{example}[thm]{Example}

\renewcommand{\leq}{\leqslant}
\renewcommand{\geq}{\geqslant}

\newcommand{\Mod}{\operatorname{Mod}}
\newcommand{\DD}{\omega_Z}
\newcommand{\Str}{\operatorname{\mathfrak{Str}}}
\newcommand{\Strw}{\operatorname{\mathfrak{Str}_{\omega}}}
\newcommand{\Ord}{\operatorname{\mathfrak{Str}_{\omega}^{\blob}}}
\newcommand{\Aw}{\mathcal B}
\newcommand{\End}{\operatorname{End}}
\newcommand{\K}{{{K}}}

\newcommand{\D}{\operatorname{\mathcal{D}}}
\newcommand{\Coh}{\operatorname{Coh}}
\newcommand{\Aut}{\operatorname{Aut}}
\newcommand{\isom}{\cong}
\renewcommand{\L}{\operatorname{L}}
\newcommand{\tensor}{\otimes}
\newcommand{\PP}{\operatorname{\mathbb P}}
\newcommand{\C}{\mathbb C}
\newcommand{\T}{\mathcal T}
\newcommand{\RR}{\mathbb{R}}

\newcommand{\blob}{{\scriptscriptstyle\bullet}}
\newcommand{\F}{\mathcal F}
\newcommand{\Fw}{{\mathcal F}_{\omega}}
\newcommand{\mat}[4]{\left( \begin{array}{cc} #1 & #2 \\ #3 & #4
\end{array} \right)}
\newcommand{\PSL}{\operatorname{PSL}}
\newcommand{\Mar}{\operatorname{\mathfrak {Mar}}}
\newcommand{\Z}{\mathbb Z}
\newcommand{\A}{\mathcal A}

\newcommand{\OO}{\mathcal O}
\newcommand{\Dw}{\operatorname{\mathcal{D}_{\omega}}}
\newcommand{\into}{\hookrightarrow}
\newcommand{\id}{\operatorname{id}}
\newcommand{\Ext}{\operatorname{Ext}}
\newcommand{\Hom}{\operatorname{Hom}}

\newcommand{\RHom}{\operatorname{Hom}}

\newcommand{\lRa}[1]{\xrightarrow{\ #1\ }}
\newcommand{\lra}{\longrightarrow}
\newcommand{\R}{\operatorname{R}}
\newcommand{\NilMod}{\operatorname{\Mod}_0}

\begin{document}

\begin{abstract}
Let $Z$ be a Fano varity satisfying the condition that the rank of the Grothendieck group of $Z$ is one more than the dimension of $Z$.
Let $\omega_Z$ denote the total space of the canonical
line bundle of $Z$,
considered as a non-compact Calabi-Yau variety. We use the theory of
exceptional collections to describe t-structures on the derived category of
coherent sheaves on $\omega_Z$. The combinatorics of these t-structures is
determined by a
natural action of an affine braid group, closely related to the well-known
action of the Artin braid group on the set of exceptional collections on
$Z$.
\end{abstract}

\maketitle

\section{Introduction}

Let $Z$ be a smooth projective Fano variety, and denote by $\omega_Z$ the total
space of its canonical bundle, which we shall think of as a non-compact
Calabi-Yau variety. The aim of this paper is to use exceptional collections
of sheaves on $Z$ to study certain sets of t-structures in the derived categories of coherent
sheaves on $Z$ and $\omega_Z$. We shall describe the combinatorics of these t-structures by introducing graphs,
whose vertices are the t-structures, and whose edges correspond to the operation
of tilting a t-structure with respect to a simple object in its heart.

It
turns out that the structure of the resulting graphs can be described using natural actions of braid groups.
The appearance of braid groups in this context is perhaps not too surprising
given the well-known action of the Artin braid group on sets of exceptional
collections discovered by Bondal \cite{B} and  Gorodentsev and Rudakov \cite{Go,GR}. In
fact Section 3 of this paper, which deals with t-structures in the derived
category of $Z$, consists of a rephrasing of part of the theory of exceptional
collections and mutations
developed by the Rudakov seminar \cite{Ru} in the language
of t-structures and tilting. Much of this story was presumably known to the
participants of this seminar.

In Section 4 we consider t-structures on the derived
category of coherent sheaves on $\omega_Z$. Our results will be used in \cite{B2} in the case $Z=\PP^2$ to describe an open subset of the space of stability conditions \cite{B} on $\omega_{\PP^2}$. Another motivation for
studying this problem is that the graphs of t-structures we construct
bear a
close resemblance to certain graphs of quiver gauge theories constructed
by the physicists Feng, Hanany, He and Iqbal \cite{FHHI}.
The edges of the physicists' graphs come from an operation
which they call Seiberg duality. We hope that studying the relationship between the
physicists' computations and the homological algebra described here will lead
to some useful insights.

Throughout we shall assume that the variety $Z$ has a full
exceptional collection  and satisfies
\begin{equation*}
\tag{$\dagger$}
\dim K(Z)\tensor \C=1+\dim Z.
\end{equation*}
Examples of such varieties include projective spaces, odd-dimensional quadrics \cite{K} and
certain Fano threefolds \cite{O}. In fact our main interest is in the case $Z=\PP^2$. Other cases not satisfying $(\dagger)$, such as
$Z=\PP^1\times\PP^1$, are more interesting and difficult, but not so
well  understood at present (see however \cite{FHH} and \cite{Ru2}).

To understand the technical significance of the assumption $(\dagger)$, recall that the class of strong exceptional collections is not closed under mutations.
On the other hand, Bondal and Polishchuk \cite{BP} introduced a class of strong exceptional
collections (see Section 3.1 for the definition),
closed under mutations,  which they referred to as geometric collections, and showed that these collections exist only on varieties satisfying $(\dagger)$. They also showed that any full exceptional collection consisting entirely of sheaves
 on such a variety is automatically geometric. We shall work with full, geometric collections throughout, but we prefer to call them \emph{simple collections}, since there is nothing particularly ungeometric about collections such as $(\OO,\OO(1,0),\OO(0,1),\OO(1,1))$ on $\PP^1\times\PP^1$ which do not satisfy Bondal and Polishchuk's conditions.

\subsection{}
Let $\D=\D^b(\Coh Z)$ denote the bounded derived category of coherent sheaves
on $Z$. Rickard's general theory of derived Morita equivalence \cite{Ri} shows that any
full, strong, exceptional
collection $(E_0,\cdots,E_{n-1})$ in $\D$ gives rise to an
equivalence of categories
\[\RHom_{\D}^{\blob}\big(\bigoplus_{i=0}^{n-1} E_i, -\big)\colon \D \lra \D^b(\Mod A),\]
where $\Mod A$ is the category of finite-dimensional right modules for the
algebra
\[A=\End_{\D} \big(\bigoplus_{i=0}^{n-1} E_i\big).\]
As explained by Bondal \cite{Bo}, the finite-dimensional algebra $A$ can be described as the path algebra of a
quiver with relations with vertices $\{0,1,\cdots,n-1\}$.
We shall always assume that the collection $(E_0,\cdots,E_{n-1})$
is a simple collection; the quiver then takes
the form
\[\xymatrix{
  {\bullet} \ar[r]^{d_1} & {\bullet} \ar[r]^{d_2}  &\bullet \ar[r] &\cdots
  \ar[r] & {\bullet} \ar[r]^{d_{n-1}}
  & {\bullet}  }
\]
with $d_i=\dim \Hom_{\D}(E_{i-1},E_i)$
arrows connecting vertex $i-1$ to vertex $i$.

Pulling back the standard t-structure on $\D(\Mod A)$ gives a t-structure on
$\D$ whose heart $\A\subset \D$ is an abelian category
equivalent to $\Mod A$. We call the subcategories
$\A\subset \D$ obtained from simple collections in this
way exceptional.
Any exceptional subcategory is of finite length and has
$n$ simple objects $S_0,\cdots ,S_{n-1}$ corresponding to the vertices of
the quiver. These simple objects have a canonical ordering coming from the
ordering of the exceptional objects $E_i$, or equivalently from the ordering
of the vertices of the quiver.

Each simple object $S_i$ defines a torsion pair in $\A$ whose torsion part
consists of direct sums of copies of $S_i$. Performing an abstract tilt in the sense of Happel, Reiten and Smal{\o} \cite{HRS}
leads to a new abelian subcategory $\L_{S_i} \A\subset \D$ which we refer to as
the left tilt of $\A$ at the simple $S_i$. It turns out that, providing $i>0$,
the category $\L_{S_i}\A\subset \D$ is also exceptional, and in fact
corresponds to a simple collection in $\D$ obtained from the
original one by a mutation. In contrast, the subcategory $\L_{S_0} \A$ has
rather strange properties in general (see Example \ref{strange}).

The fact that mutations of exceptional collections give rise to an
action of the Artin braid group now translates as

\medskip

\noindent {\bf Theorem \ref{one} }{\it 
The Artin braid group $A_n$ acts on the set of exceptional
subcategories of $\D$. For each integer $1\leq i<n-1$ the generator $\sigma_i$
acts by tilting a subcategory at its $i$th simple object.}

\medskip

It is convenient to introduce a graph $\Str(Z)$ whose vertices are exceptional
subcategories of $\D$, and in which two vertices are linked by an edge if the
corresponding abelian subcategories are related by a tilt at a simple object. In the case $Z=\PP^2$ we shall show that the action of Theorem \ref{one} is free.
It follows that  each connected component of $\Str(\PP^2)$ is the Cayley
graph of the standard system of generators of the group $A_n$.

\subsection{}
Consider now the category $\D^b(\Coh \omega_Z)$. Any simple exceptional
collection $(E_0,\cdots,E_{n-1})$ in $\D$ determines an equivalence
\[\RHom_{\DD}^{\blob}\big(\bigoplus_{i=0}^{n-1} \pi^* E_i, -\big)\colon \D^b(\Coh \omega_Z) \lra \D^b(\Mod B),\]
where $\pi\colon \omega_Z \to Z$ is the projection, and $\Mod B$ is the category of finitely generated right modules for the
algebra
\[B=\End_{\DD} \big(\bigoplus_{i=0}^{n-1} \pi^* E_i\big).\]
Note that the algebra $B$ is infinite-dimensional. Nonetheless $B$
can again be described as the path algebra of a quiver with relations with vertices $\{0,1,\cdots,n-1\}$. This
time the quiver is of the form
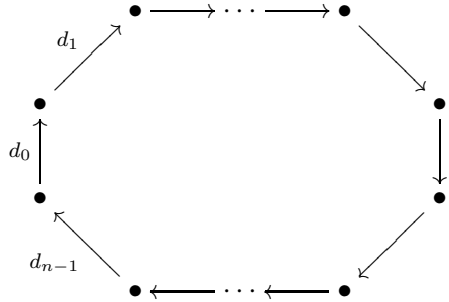
\[\xymatrix{ &  \bullet \ar[r] &\cdots \ar[r] & \bullet \ar[dr] \\
\bullet \ar[ur]^{d_1} &&&& \bullet \ar[d]
\\ \bullet \ar[u]^{d_{0}}
 &&&& \bullet \ar[dl] \\ & \bullet \ar[ul]^{d_{n-1}} &\cdots \ar[l] &\bullet \ar[l]}\]
with $d_i$ arrows from vertex $i-1$ to vertex $i$ for $1\leq i\leq n-1$ as before, and
\[d_0=\dim\Hom_{\D}(E_{n-1}\tensor \omega_Z,E_0)\]
arrows connecting vertex $n-1$ to vertex $0$.

Consider the full subcategory $\Dw\subset\D^b(\Coh \omega_Z)$ consisting of objects
supported on the zero section $Z\subset \omega_Z$. The above equivalence
determines a t-structure on $\Dw$ whose heart is an abelian subcategory
$\Aw\subset\Dw$ equivalent to the category of nilpotent
representations of the algebra $B$. Abelian subcategories $\Aw\subset\Dw$
obtained in this way will be again be called exceptional.
Any exceptional subcategory of $\Dw$ is of
finite length and has $n$ simple objects $S_0,\cdots,S_{n-1}$ corresponding
to the vertices of the quiver. These simple objects have a canonical ordering
coming from the ordering of the exceptional objects $(E_0,\cdots, E_{n-1})$, and
for $1\leq i< n-1$, the abelian subcategory $\L_{S_i} \Aw\subset \Dw$ is
also exceptional, and corresponds to a simple collection
 in $\D$ obtained from the original one by a mutation.

The key new feature of the Calabi-Yau situation concerns the subcategory $\L_{S_0} \Aw$. The simple objects $S_i$ of an exceptional subcategory $\Aw\subset\Dw$ are
 spherical objects. It follows from work of Seidel and Thomas \cite{ST}
that there are associated autoequivalences $\Phi_{S_i}\in \Aut \Dw$, and we shall show that the category $\L_{S_0} \Aw\subset\Dw$ is the image of an
exceptional subcategory of $\Dw$ under the autoequivalence $\Phi_{S_0}$.

A
subcategory $\Aw\subset\Dw$ will be called quivery if there is an
autoequivalence $\Phi\in \Aut\Dw$ such that the subcategory
$\Phi(\Aw)\subset \Dw$ is exceptional. Thus, quivery subcategories of $\Dw$
are finite length abelian categories, and from what was said above,
they remain quivery under the operation of
tilting at a simple object.
A slightly subtle point is that the simple objects $S_0,\cdots,S_{n-1}$ of a
quivery subcategory $\Aw\subset\Dw$ have no canonical ordering, only
a cyclic ordering coming from the arrows in the corresponding quiver.
Let us define an
ordered quivery subcategory to be a quivery subcategory
$\Aw\subset \Dw$ together with an ordering of its $n$ simple objects
$(S_0,\cdots,S_{n-1})$ compatible with the canonical cyclic ordering.

The combinatorics of the set of quivery subcategories of $\Dw$ is controlled not by the Artin braid group $A_n$, but by a group $B_n$ which is a quotient of the annular braid group $CB_n$,
or alternatively, a semidirect product of the affine braid group $\tilde{A}_{n-1}$ by the cyclic group $\Z_n$. The reader is referred to Section 2.1 for the precise definitions of these groups. 
 
\medskip

\noindent {\bf Theorem \ref{two}. }{\it 
There is an action of the group $B_n$ on
the set of ordered quivery subcategories of $\Dw$. For each integer
$0\leq i\leq n-1$ the element $\tau_i$ acts on the underlying abelian subcategories by
tilting at the $i$th simple object.}

\medskip

Introduce a graph $\Strw(Z)$ whose vertices are the quivery
subcategories of $\Dw$, and in which two vertices are joined by an
edge if the corresponding subcategories are related by a tilt at a
simple object. In the case $Z=\PP^2$ we shall show that the action of
Therorem \ref{two} is free, and it follows that each connected component of the graph
$\Strw(\PP^2)$ is the Cayley graph for the standard system of generators $\tau_0,\cdots, \tau_{n-1}$ of the affine braid group
$\tilde{A}_{n-1}$.

\subsection*{Acknowledgements}
This paper has benefitted from conversations I had with many mathematicians; I'd particularly like to thank Phil Boalch, Alastair King and Michel van den Bergh.


\section{Preliminaries: Braid groups and tilting}

This section consists of various basic facts and definitions we shall need;
we include
the material here for the reader's convenience, and to fix notation.

\subsection{Braid groups}
Given a topological space $M$, define the $n$-point configuration space
\[ C_n(M)=\big\{(m_0,\cdots,m_{n-1})\in M^n: i\neq j\implies m_i\neq m_j\big\}.\]
The symmetric group $\Sigma_n$ acts freely on $C_n(M)$ permuting the points.

The standard $n$-string \emph{Artin braid group} $A_n$ is defined to be the fundamental group of the space $C_n(\C)/\Sigma_n$.
As is well-known (see for example \cite{Bi}), it is generated by elements $\sigma_1,\cdots,\sigma_{n-1}$
subject to the relations \begin{eqnarray*}
 \sigma_i\sigma_{i+1}\sigma_i
&=&\sigma_{i+1}\sigma_i\sigma_{i+1} \text{ for } 1\leq i< n-1, \\
\sigma_i\sigma_j&=&\sigma_j\sigma_i \text{ for } j-i\neq\pm 1.
\end{eqnarray*}
The centre of $A_n$ is generated by the element
\[\gamma=(\sigma_1\cdots \sigma_{n-1})^n=(\sigma_{n-1}\cdots \sigma_1)^n.\]
To visualize elements of the group $A_n$  one can project points in $C_n(\C)$ to a far away line in $\C$ to obtain a set of $n$ points in $\RR$; a loop in the configuration space can
then be thought of as a braid on $n$ strings. The elementary generators $\sigma_i$ correspond to the $i$th string passing under the $(i-1)$st. 

We shall need the following easy result later.

\begin{lemma} \label{conj} The element
\[ \delta= (\sigma_1\cdots
\sigma_{n-1})(\sigma_1\cdots \sigma_{n-2})\cdots
(\sigma_1\sigma_2)\sigma_1\in A_n\] has the property that
$\delta^{-1}\sigma_i\delta=\sigma_{n-i}$ for $1\leq i\leq n-1$.
\end{lemma}

\begin{pf}
For $1\leq j\leq n-1$ set $\beta_j=\sigma_1\cdots \sigma_j$. We are required
to prove that
\[\sigma_i\,\beta_{n-1}\beta_{n-2}\cdots\beta_1=
\beta_{n-1}\beta_{n-2}\cdots\beta_1\,\sigma_{n-i}.\] First suppose $i>1$. By
induction on $n$ we can assume that
\[\sigma_{i-1}\,\beta_{n-2}\cdots\beta_1=\beta_{n-2}\cdots\beta_1\,\sigma_{n-i}.\]
Multiplying both sides by $\beta_{n-1}$ and noting that for $1<i\leq n-1$ we
have $\beta_{n-1}\sigma_{i-1}=\sigma_i\,\beta_{n-1}$ gives the result. To
prove the result when $i=1$ note first that $\sigma_{n-1}$ commutes with
$\beta_j$ if $j\leq n-3$. Thus we are reduced to proving
\[\sigma_1\beta_{n-1}\beta_{n-2}=\beta_{n-1}\beta_{n-1}.\]
This follows by repeatedly applying the relation
$\sigma_i\,\beta_{n-1}=\beta_{n-1}\sigma_{i-1}$.
\end{pf}

The $n$-string ($n\geq 2$) \emph{annular braid group} is defined to be the fundamental group of the space $C_n(\C^*)/\Sigma_n$.
It is generated by elements $\tau_i$ indexed by the cyclic group
$\Z_n$, together with a single element $r$, subject to the relations
\begin{eqnarray*}
r\tau_i r^{-1}&=&\tau_{i+1} \text{ for all }i\in \Z_n, \\
\tau_i\tau_{i+1}\tau_i&=&\tau_{i+1}\tau_i\tau_{i+1}\text{ for all }i\in
\Z_n, \\
\tau_i\tau_j&=&\tau_j\tau_i\text{ for }j-i\neq \pm 1.\\
\end{eqnarray*}
For a proof of the validity of this presentation see \cite{KP}. Of more interest to us will be the quotient group
\[B_n=CB_n/\langle r^n\rangle.\]
The subgroup of $B_n$ (or $CB_n$) generated by the elements $\tau_0, \cdots
,\tau_{n-1}$ is an affine braid group; we denote it $\tilde{A}_{n-1}$.

To visualize elements of these groups one can project points in $C_n(\C^*)$ out from the origin onto a large circle to obtain $n$ points in $S^1$; a loop in the configuration space can then be thought of as a braid of $n$ strings lying on
the surface of a cylinder. The element $\tau_i$ corresponds to the $i$th string passing under the $(i-1)$st; the element $r$ corresponds to the twist which for each $i$ takes point $i$ to point $i+1$.

\begin{prop}
\label{exactsequence} There is a short exact sequence
\[ 1\lra F_n\lra B_n\lRa{h} A_n/\langle\gamma\rangle\lra 1,\]
where $F_n$ is the free group on $n$ generators. The homomorphism $h$ is defined by
\[h(r)=\sigma_1\cdots\sigma_{n-1}\text{ and }h(\tau_i)=\sigma_i\text{ for }1\leq i\leq n-1 ,\]
and its kernel is freely generated by the elements
\[\alpha_i=r^i(\tau_1\cdots \tau_{n-1})r^{-(i+1)}\qquad 0\leq i\leq n-1.\]
 \end{prop}

\begin{pf}
We give two proofs, one geometric and the other algebraic. In geometric terms, note that the space $C_n(\C^*)/\Sigma_n$ is homotopic to $C_{n+1}(\C)/\Sigma_n$
where $\Sigma_n\subset \Sigma_{n+1}$ is the subgroup fixing $n\in\{0,1,\cdots, n\}$.
Forgetting the last point gives a fibration
\[C_{n+1}(\C)/\Sigma_n\lra C_n(\C)/\Sigma_n\]
whose fibre is $\C\setminus\{m_0,\cdots ,m_{n-1}\}$. This gives an exact sequence
\[1\lra F_n \lra CB_n\lRa{h} A_n\lra 1.\]
Drawing suitable pictures it is easy enough to see that $h$ acts on generators as claimed in the statement, and that the elements $\alpha_i$ correspond to loops in the fibre which freely generate the fundamental group of  $\C\setminus\{m_0,\cdots ,m_{n-1}\}$.
Since $h(r^n)=\gamma$ the result follows by taking quotients.
 
To see the result using just the presentation of $B_n$ we follow an argument of Chow \cite{C}. It is easy to check that the formula in the statement defines a homomorphism $h\colon CB_n \to A_n$, and that the elements $\alpha_i$ lie in its kernel and generate a normal subgroup $K\subset CB_n$. Furthermore $h$ has a section $A_n \to CB_n$ sending $\sigma_i$ to $\tau_i$ for $1\leq i\leq n-1$, and the induced homomorphism $A_n\to CB_n/K$ is surjective because in $CB_n/K$ one has $r=\tau_1\cdots \tau_{n-1}$. It follows that $K$ is the kernel of $h$.

The only non-trivial part is to show that $K\subset CB_n$ is freely generated by the elements $\alpha_i$. To see this, one needs to exhibit a representation of $CB_n$ in which they act freely. Let $F_n$ be the free group on generators $x_i$ indexed by $i\in \Z_n$, and define an action of  $CB_n$ on $F_n$ by automorphisms using the formulae
$r(x_i)=x_{i+1}$ and
\[\tau_i(x_i)=x_{i+1},\quad \tau_i(x_{i+1})=x_{i+1}^{-1}x_i x_{i+1},\quad \tau_i(x_j)=x_j \text{ for }j\notin \{i,i+1\}.\]
Then the element $\alpha_i$ acts by sending each $x_j$ to $x_i x_j x_i^{-1}$
and it follows that the $\alpha_i$ generate the free group of inner automorphisms of $F_n$.
\end{pf}

\subsection{T-structures and tilting}
The reader is assumed to be familiar with the concept of a t-structure
\cite{BBD,GM}. The following easy result is a good exercise.

\begin{lemma}
\label{ll} A bounded t-structure is determined by its heart. Moreover, if
$\A\subset\D$ is a full additive subcategory of a triangulated category $\D$,
then $\A$ is the heart of a bounded t-structure on $\D$ if and only if the
following two conditions hold:
\begin{itemize}
\item[(a)]if $A$ and $B$ are objects of
$\A$ then $\Hom_{\D}(A,B[k])=0$ for $k<0$,

\item[(b)]
for every nonzero object $E\in\D$ there are integers $m<n$ and a collection
of triangles
\[
\xymatrix@C=.2em{ 0_{\ } \ar@{=}[r] & E_{m} \ar[rrrr] &&&& E_{m+1} \ar[rrrr]
\ar[dll] &&&& E_{m+2} \ar[rr] \ar[dll] && \ldots \ar[rr] && E_{n-1}
\ar[rrrr] &&&& E_n \ar[dll] \ar@{=}[r] &  E_{\ } \\
&&& A_{m+1} \ar@{-->}[ull] &&&& A_{m+2} \ar@{-->}[ull] &&&&&&&& A_n
\ar@{-->}[ull] }
\]
with $A_i[i]\in\A$ for all $i$.\qed
\end{itemize}
\end{lemma}

It follows from the definition that the heart of a bounded t-structure is an abelian category \cite{BBD}.
In analogy with the standard t-structure on the derived category of an
abelian category, the objects $A_i[i]\in \A$ are called the cohomology
objects of $A$ in the given t-structure, and denoted
$H^i(E)$.

Note that the group $\Aut \D$ of exact autoequivalences of $\D$ acts on the set of
bounded t-structures: if $\A\subset\D$ is the heart of a bounded t-structure
and $\Phi\in \Aut \D$, then $\Phi(\A)\subset \D$
is also the heart of a bounded t-structure.

A very useful way to construct t-structures is provided by the method of
tilting. This was first introduced in this level of generality
by Happel, Reiten and Smal{\o} \cite{HRS}, but the name and the basic idea go back to a paper of Brenner and Butler \cite{BB}.

\begin{defn}
\label{tors} A torsion pair in an abelian category $\A$ is a pair of full
subcategories $(\T,\F)$ of $\A$ which satisfy $\Hom_{\A}(T,F)=0$ for $T\in
\T$ and $F\in\F$, and such that every object $E\in\A$ fits into a  short
exact sequence
\[0\lra T\lra E\lra F\lra 0\] for some pair of objects $T\in\T$ and
$F\in  \F$.
\end{defn}

The objects of $\T$ and $\F$ are called torsion and torsion-free. The proof
of the following result \cite[Proposition 2.1]{HRS} is pretty-much immediate from Lemma \ref{ll}.

\begin{prop}(Happel, Reiten, Smal{\o})
Suppose $\A$ is the heart of a bounded t-structure on a triangulated category
$\D$. Given an object $E\in\D$ let $H^i(E)\in\A$ denote the $i$th cohomology
object of $E$ with respect to this t-structure. Suppose $(\T,\F)$ is a
torsion pair in $\A$. Then the full subcategory
\[\A^{\sharp}=\big\{E\in \D:H^i(E)=0\text{ for }i\notin\{-1,0\},
H^{-1}(E)\in\F\text{ and } H^{0}(E)\in\T\big\}\] is the heart of a bounded
t-structure on $\D$.\qed
\end{prop}

In the situation of the Lemma one says that the the subcategory $\A^{\sharp}$ is
obtained from the subcategory $\A$ by tilting with respect to the
torsion pair $(\T,\F)$. In
fact one could equally well consider $\A^{\sharp}[-1]$ to be the tilted
subcategory; we shall be more precise about this where necessary. Note that the
pair $(\F[1],\T)$ is a torsion pair in $\A^{\sharp}$ and that tilting with
respect to this pair gives back the original subcategory $\A$ with a shift.

\medskip

Now suppose $\A\subset \D$ is the heart of a bounded t-structure and is a
finite length abelian category. Note that the t-structure is completely determined by the set of simple objects of $\A$; indeed $\A$ is the smallest extension-closed subcategory of $\D$ containing this set of objects. Given a simple object $S\in \A$ define
$\langle S \rangle\subset \A$ to be the full subcategory consisting of objects $E\in\A$ all
of whose simple factors are isomorphic to $S$. One can either view $\langle S \rangle$ as
the torsion part of a torsion theory on $\A$, in which case the torsion-free
part is
\[\F=\{E\in \A:\Hom_{\A}(S,E)=0\},\]
or as the torsion-free part, in which case the torsion part is
\[\T=\{E\in\A:\Hom_{\A}(E,S)=0\}.\]
The corresponding tilted subcategories are
\begin{eqnarray*}
\L_S \A &=& \{E\in\D:H^i(E)=0\text{ for
}i\notin\{0,1\},H^{0}(E)\in\F\text{ and }H^1(E)\in\langle S \rangle\}  \\
\R_S \A &=& \{E\in \D:H^i(E)=0\text{ for
}i\notin\{-1,0\},H^{-1}(E)\in\langle S \rangle\text{ and }H^0(E)\in\T\}.
\end{eqnarray*}
We define these subcategories of $\D$ to be the left  and right tilts of the subcategory
$\A$ at the simple $S$ respectively. It is easy to see that $S[-1]$ is a simple object of
$\L_S \A$, and that if this category is finite length, then $\R_{S[-1]} \L_S
\A=\A$. Similarly, if $\R_S \A$ is finite length $\L_{S[1]}\R_S \A=\A$.

\smallskip

The following obvious result will often be useful.

\begin{lemma}
\label{useful} The operation of tilting commutes with the action of the group
of autoequivalences on the set of t-structures. Take an autoequivalence $\Phi\in \Aut \D$. If $\A\subset\D$ is the
heart of a bounded t-structure on $\D$ and has finite length and $S\in \A$ is simple,
then $\Phi(\A)\subset \D$ is the heart of a bounded t-structure on $\D$ and has finite length, $\Phi(S)$ is a simple object of $\Phi(\A)$, and
\[\L_{\Phi(S)} \Phi(\A)=\Phi(\L_S \A).\]
\end{lemma}

\begin{pf}
This is a straightforward application of the definitions.
\end{pf}


\section{Exceptional collections and t-structures on $\D$}

Throughout this section $Z$ will be a smooth projective Fano variety
and $\D$ will be its bounded derived category of coherent sheaves.
We shall assume throughout that $Z$ satisfies the condition
\[\tag{$\dagger$} \dim \K(Z)\tensor \C=1+\dim Z.\]
Although this is not necessary everywhere, some of the definitions would need to be modified for more general cases, and it is not clear exactly how this should be done.

\subsection{Exceptional collections and mutations}
We start by recalling some of the theory of exceptional collections developed by Bondal, Gorodentsev, Polishchuk, Rudakov and others. For more information
and proofs of some of the following facts the
reader is referred to the original papers \cite{Bo,BP,Go,GR,Ru}.

An object $E\in \D$ is said to be \emph{exceptional} if
\[\Hom_{\D}^k(E,E)=\biggl\{\begin{array}{ll} \C & \text{ if }k=0, \\ 0 &\text{
otherwise.}\end{array}\biggr.\] An \emph{exceptional collection} in $\D$
(or on $Z$) of
length $n$ is a sequence of exceptional objects $(E_0,\cdots ,E_{n-1})$ of
$\D$ such that
\[ n-1\geq i>j \geq 0\implies \Hom^k_{\D}(E_i,E_j)=0\text { for all }k\in\Z.\]
The exceptional collection $(E_0,\cdots ,E_{n-1})$ in $\D$ is \emph{full} if for any
$E\in \D$
\[\Hom^k_{\D}(E_i,E)=0\text{ for all }0\leq i\leq n-1\text{ and all
}k\in\Z\implies E\isom 0.\]
An exceptional collection $(E_0,\cdots,E_{n-1})$ is \emph{strong} if
for all $0\leq i,j\leq n-1$ one has\[\Hom_{\D}^k(E_i,E_j)=0\text{ for }k\neq 0.\]
As we shall see in the next subsection, strong exceptional collections
define equivalences of $\D$ with derived categories of module categories.
Pulling back the standard t-structure allows us to define new t-structures on $\D$. Thus
if we are interested in t-structures on $\D$ exceptional collections are not enough: we need strong collections.

\smallskip

Given two objects $E$ and $F$ of $\D$, define a third object $\L_E F$ of $\D$ (up to isomorphism) by the triangle
\[\L_E F \lra \Hom_{\D}^{\blob}(E,F)\tensor E\lRa{ev} F,\]
where $ev$ denotes the canonical evaluation map.
It is easy to see that if $(E,F)$ is an exceptional collection then so is $(\L_E F,E)$. 
The object $\L_E F$ is called
the \emph{left mutation of $F$ through
$E$}. Mutations of this form define a
braid group action on exceptional collections \cite{Bo,Go,GR}.

\begin{thm}(Bondal, Gorodentsev, Rudakov)
\label{exceptionalaction} The braid group $A_n$ acts on the set of
exceptional collections of length $n$ in $\D$ by mutations. For $1\leq i\leq
n-1$, the generating element $\sigma_i$ acts by
\[\sigma_i(E_0,\cdots ,E_{n-1}) = (E_0,\cdots, E_{i-2}, \L_{E_{i-1}}
E_i,E_{i-1}, E_{i+1},\cdots ,E_{n-1}).\qed\]
\end{thm}

Strong exceptional collections do not remain strong under mutations in general. A good example is the
strong collection $(\OO,\OO(1,0),\OO(0,1),\OO(1,1))$ on $\PP^1\times \PP^1$ which mutates to give the non-strong collection $(\OO,\OO(0,1)[-1],\OO(1,0),\OO(1,1))$.

\smallskip

A \emph{helix} in $\D$ is an infinite sequence of objects $(E_i)_{i\in \Z}$
such that for each $i\in \Z$ the corresponding \emph{thread} $(E_i, \cdots ,
E_{i+n-1})$ is a full exceptional collection in $\D$, and
the relation
\[ (\sigma_1\cdots \sigma_{n-1}) (E_{i+1}, \cdots ,E_{i+n})=(E_i, \cdots,
E_{i+n-1})\] is satisfied. Clearly a helix $(E_i)_{i\in \Z}$ is uniquely determined by the full exceptional collection $(E_0,\cdots, E_{n-1})$; we say that the helix is \emph{generated} by  $(E_0,\cdots, E_{n-1})$. Bondal \cite[Theorem 4.2]{Bo} showed that any helix $(E_i)_{i\in \Z}$
satisfies
\begin{equation}
\label{bondal}
E_{i-n}\isom E_i\tensor\omega_Z\text{ for all }i\in \Z.
\end{equation}
These definitions certainly need to be modified for varieties
$Z$ not satisfying $(\dagger)$, but it is not clear exactly how this should be done.

We shall call a helix $(E_i)_{i\in \Z}$ in $\D$ \emph{simple} if for all $i\leq j$
one has
\[\Hom_{\D}^k (E_i, E_j) = 0\text{ unless } k=0.\]
Such helices were called \emph{geometric} by Bondal and Polishchuk. 
An exceptional collection $(E_0,\cdots,E_{n-1})$ will be called \emph{simple} if it is a full collection which generates a simple helix.
Equivalently this means that the collection is full, and for any integers $0\leq i,
j\leq n-1$ and any $p\leq 0$
\[\Hom_{\D}^k(E_i,E_j\tensor\omega_Z^p)=0 \text{ unless }k=0.\]
In particular, any simple collection is strong.
Bondal and Polishchuk showed that any full exceptional collection
of sheaves on a variety satisfying $(\dagger)$ is automatically simple
\cite[Proposition 3.3]{BP}.

The importance of simple collections is the following result \cite[Theorem 2.3]{BP}.

\begin{thm}(Bondal, Polishchuk)
\label{BondalPolishchuk}
Any mutation of a simple collection is again simple.
\qed
\end{thm}

The motivating example for all this theory is the
sequence of line bundles \[(\OO,\OO(1), \cdots ,\OO(n-1))\] on $\PP^{n-1}$, which is a simple
collection of length $n$. The fact that it is full is the
essential content of Beilinson's theorem \cite{Be}.
The helix generated by this collection  is just $(\OO(i))_{i\in \Z}$.

\subsection{The homomorphism algebra}
Let $(E_0,\cdots ,E_{n-1})$ be a full, strong exceptional collection in $\D$.
The general theory of derived Morita equivalence \cite{Ri} shows that the functor
\[ \F=\RHom_{\D}^{\blob}\big(\bigoplus_{i=0}^{n-1} E_i, -\big)\colon \D \lra \D(\Mod A)\]
is an equivalence, where $\Mod A$ is the category of finite-dimensional right
modules for the algebra
\[A=\End_{\D} \big(\bigoplus_{i=0}^{n-1} E_i\big).\]
This algebra is called the \emph{homomorphism algebra} of the collection $(E_0,\cdots,E_{n-1})$. 
Note that $A$ is finite-dimensional and has a natural grading
\[A=\bigoplus_{k=0}^{n-1} A_k=\bigoplus_{k=0}^{n-1} \bigoplus_{j-i=k} \Hom_{\D}(E_i,E_j).\]
The degree zero part has a basis consisting of the idempotents
\[e_i=\id_{E_i}\in \End_{\D}(E_i),\]
and there are corresponding simple right-modules $T_0,\cdots,T_{n-1}$ defined by
\[\dim_{\C} (T_j e_i)=\delta_{ij}.\]
It is easy to check that all simple modules are of this form.

\begin{prop}(Bondal)
\label{baah}
Let $(E_0,\cdots, E_{n-1})$ be a full, strong exceptional collection in $\D$, and define a new collection by
\[(F_0, \cdots , F_{n-1})=\delta (E_0, \cdots E_{n-1}),\]
where $\delta\in A_n$ is the element defined in Lemma \ref{conj}.
Then these two collections are dual, in the sense that
\[ \Hom^k_{\D}(E_i, F_{n-1-j}[j])=
\biggl\{\begin{array}{ll} \C & \text{ if }i=j\text{ and } k=0,\\ 0 &\text{
otherwise.}\end{array}\biggr.\]
The objects $F_i$ are unique with this property.
\end{prop}

\begin{pf}
This is basically Lemma 5.6 of \cite{Bo}. Just note that in
Bondal's notation
\[\delta(E_0,\cdots,E_{n-1})=(\L_{n-1} E_{n-1},\cdots,\L_1 E_1, E_0),\]
where for $1\leq i\leq n-1$ the object $\L_i E_{i}$ is defined to be $\L_{E_0}\L_{E_1}\cdots\L_{E_{i-1}} E_i$.
\end{pf}

Under the equivalence $\F$, the object $E_i\in \D$ is mapped to the
projective module $e_i A$ corresponding to the vertex $i$.
Lemma \ref{baah} shows that the object \[S_j=F_{n-1-j}[j]\]
is mapped to the simple module
$T_j$. Note also that Lemma \ref{conj} shows that mutations of the collections $(E_0,
\cdots, E_{n-1})$ and $(F_0,\cdots,F_{n-1})$ correspond to each other.

As an example, take the collection 
$(\OO,\OO(1),\cdots,\OO(n-1))$ in $\D(\PP^{n-1})$. The dual collection, in the sense of Lemma \ref{baah}, is \[(\Omega^{n-1}(n-1),\cdots,
\Omega^{1}(1),\OO),\] where $\Omega^i=\bigwedge^i T^*$ is the sheaf of holomorphic $i$-forms on $\PP^{n-1}$. This can be checked directly by computing the cohomology groups of Proposition \ref{baah}.

\begin{prop}(Bondal, Polishchuk)
\label{generated}
 Let $(E_0,\cdots, E_{n-1})$ be a simple collection in $\D$ and let $A$ be the corresponding homomorphism algebra with its natural grading. Then
 $A$ is generated over $A_0$ by $A_1$ and is
 Koszul.
\end{prop}

\begin{pf}
For the first statement it is enough to show that for $0\leq i<j\leq
n-1$, the natural map
\[\Hom_{\D}(E_i,E_{j-1})\tensor \Hom_{\D}(E_{j-1},E_j)\lra \Hom_{\D}(E_i,E_j)\]
is surjective. Thus it is enough to show that
\[\Hom_{\D}^1(E_i,\L_{E_{j-1}} E_j)=0.\]
This statement follows from the fact that the collection $\sigma_j(E_0,\cdots,E_{n-1})$ is strong, 
which in turn follows from Theorem \ref{BondalPolishchuk}.

The condition that $A$ is Koszul is equivalent to the statement that the Yoneda algebra
\[A^!=\End_{A}^{\blob}\big(\bigoplus_{j=0}^{n-1} T_j\big)\]
is generated in degree one. Under the equivalence $\F$ described above, the simple modules $T_j$ correspond to the objects $S_j=F_{n-1-j}[j]$. Thus $A^!$ is just the homomorphism algebra of the dual exceptional collection  $(F_0,\cdots,F_{n-1})$. By Theorem \ref{BondalPolishchuk} this collection is also simple, so the result follows.
\end{pf}

The homomorphism algebra of a simple collection can  naturally be thought of as the
path algebra of a quiver with relations. The quiver has $n$ vertices
$\{0,1,\cdots ,n-1\}$ corresponding to the idempotents $e_i$, and for each
$1\leq i\leq n-1$ has
\[d_i=\dim \Hom_{\D}(E_{i-1},E_{i})\] arrows going from vertex
$i-1$ to vertex $i$.
\[\xymatrix{
  {\bullet} \ar[r]^{d_1} & {\bullet} \ar[r]^{d_2}  &\bullet \ar[r] &\cdots
  \ar[r] & {\bullet} \ar[r]^{d_{n-1}}
  & {\bullet}  }
\]
Since the algebra is Koszul the relations are quadratic \cite{BGS}.

\subsection{Tilting and mutations}

Given a simple collection $(E_0,\cdots, E_{n-1})$ in $\D$, the
corresponding equivalence
\[ \F=\RHom_{\D}^{\blob}\big(\bigoplus_{i=0}^{n-1} E_i, -\big)\colon \D \lra \D(\Mod A)\]
allows one to pull back the standard t-structure on $\D(\Mod A)$ to give a
t-structure on $\D$ whose heart \[\A(E_0,\cdots, E_{n-1})\subset \D\]
is
equivalent to the abelian category $\Mod A$. Let us call the
subcategories of $\D$ obtained in this way \emph{exceptional}. Note that any
exceptional subcategory  is a
finite length abelian category with $n$ simples $S_0,\cdots,S_{n-1}$.
These simples have a uniquely defined ordering $(S_0,\cdots,S_{n-1})$ in
which
\begin{equation}
\label{slip}
\Hom_{\D}^k(S_i,S_j)=0 \text{ unless } i-j=k\geq 0.\
\end{equation}
Thus it is possible to talk about the $i$th simple object $S_i$ of an
exceptional subcategory.

\begin{prop}
\label{helloo}
Let $(E_0,\cdots,E_{n-1})$ be a simple collection in $\D$, and let $S_i$
denote the $i$th simple object of the exceptional
subcategory $\A(E_0,\cdots,E_{n-1})\subset \D$.  Then for each integer $1\leq i\leq n-1$ there is an
identification of subcategories of $\D$
\[\L_{S_i}\A(E_0,\cdots,E_{n-1})=\A(\operatorname{\sigma_i}(E_0,\cdots,E_{n-1})).\]
\end{prop}

\begin{pf}
Put $(E_0',\cdots,E_{n-1}')=\sigma_i(E_0,\cdots,E_{n-1})$ and set
\[\A=\A(E_0,\cdots,E_{n-1}),\qquad \A'=\A(E'_0,\cdots,E'_{n-1}).\]
Let $(S_0,\cdots,S_{n-1})$ be the simple objects of $\A$ with their canonical
ordering. The subcategory $\L_{S_i} \A$ is obtained by tilting $\A$ with respect to
the torsion theory $(\T,\F)$, where $\T$ consists of direct sums of $S_i$,
and
\[\F=\{E\in \A:\Hom_{\A}(S_i,E)=0\}.\] Note that $S_j\in \F$ for every $j\neq
i$.
It will be enough to show that $\A'\subset \L_{S_i} \A$, because if two bounded t-structures have nested hearts then they are the same. Since $\A'$ has
finite length it will be enough to show that every simple object of $\A'$ is
contained in either $\T[-1]$ or in $\F$.

Recall that if $(F_0,\cdots,F_{n-1})$ is the
dual exceptional collection to $(E_0,\cdots, E_{n-1})$ then
$S_j=F_{n-1-j}[j]$. Let $(S_0',\cdots,S_{n-1}')$ be the simple objects of $\A'$ with their
canonical ordering. By Lemma \ref{conj}, the dual collection to
$(E_0',\cdots,E_{n-1}')$ is
\[(F_0',\cdots,F_{n-1}')=\sigma_{n-i}(F_0,\cdots,F_{n-1}),\]
and $S_j'=F_{n-j-1}'[j]$. For $j\notin\{i-1,i\}$ we have
$S_j'=S_j$ so that $S_j'\in \F$. Furthermore, $S_{i-1}'=S_i[-1]$. Thus the
only thing to check is that $S_i'\in \F$.

Now
$F_{n-i-1}'=\L_{F_{n-i-1}}F_{n-i}$, and rewriting the defining triangle
\[\L_{F_{n-i-1}}F_{n-i}\lra \Hom_{\D}^{\blob}(F_{n-i-1},F_{n-i})\tensor F_{n-i-1}\lRa{ev} F_{n-i},\]
we obtain a triangle
\[ \Hom_{\D}^1(S_{i},S_{i-1})[-1]\tensor S_i \lRa{ev} S_{i-1}\lra S_i',\]
where we have used (\ref{slip}) to see that $\Hom_{\D}^{\blob}(S_i,S_{i-1})$ is concentrated in degree 1.
Rewriting this triangle again shows that $S_i'$ is a universal extension in $\A$
\[0\lra S_{i-1}\lra S_i'\lra \Ext^1_{\A}(S_i,S_{i-1})\tensor S_i\lra 0,\]
and applying the functor $\Hom_{\D}(S_i,-)$ it follows that $S_i'\in\F$.
\end{pf}

Using this Lemma the braid group action on exceptional collections described in Lemma \ref{exceptionalaction} can be translated into the following form.

\begin{thm}
\label{one} The Artin braid group $A_n$ acts on the set of exceptional
subcategories of $\D$. For each integer $1\leq i<n-1$ the generator $\sigma_i$
acts by tilting a subcategory at its $i$th simple object.\qed
\end{thm}

As a final remark in this section, suppose $\A\subset \D$ is an exceptional subcategory  of $\D$ with corresponding ordered simple objects $(S_0,\cdots,S_{n-1})$. The categories $\L_{S_0}
\A$ and $\R_{S_{n-1}} \A$ are not covered by the above results. In general
these categories are rather strange, as the following example shows.

\begin{example}
\label{strange} Consider the case $\A=\A(\OO,\OO(1))\subset \D(\PP^1)$
corresponding to the simple collection $(\OO,\OO(1))$ on $\PP^1$. The dual
collection is $(\OO(-1),\OO)$ so that the simple objects of $\A$ are
$S_0=\OO$ and $S_1=\OO(-1)[1]$. The only objects $E\in \A$ satisfying $\Hom_{\A}(S_0,E)=0$ are direct sums of copies of $S_1$. Performing a left tilt at the simple $S_0$
leads to a category $\L_{\OO} \A$ which is finite length and has two simple
objects $S'_0=\OO[-1]$ and $S'_1=\OO(-1)[1]$.
Since
\[\Ext^1_{\A'}(S'_0,S'_1)=0=\Ext^1_{\A'}(S'_1,S'_0),\]
the category $\A'$ is semisimple, and so every object in the derived category
$\D(\A')$ is a direct sum of copies of shifts of $S'_0$ and $S'_1$. In
particular, the only exceptional objects in $\D(\A')$ are shifts of $S'_0$ and
$S'_1$. It follows immediately that $D(\A')$ is not equivalent to $\D$, so that the bounded t-structure whose heart is $\A'$ is unfaithful.
\end{example}


\section{Spherical collections and t-structures on $\Dw$}

Recall our general assumption: $Z$ is a smooth projective Fano variety satisfying
\[\dim K(Z)\tensor
\C=1+\dim Z,\]and $\omega_Z$ is the canonical bundle of $Z$, which we view
both as an invertible $\OO_Z$-module, and as a quasi-projective variety with
a fibration $\pi\colon \omega_Z\to Z$. The inclusion of the zero section in
$\omega_Z$  will be denoted $s\colon Z\into \omega_Z$.
Define
\[\Dw\subset \D^b(\Coh \omega_Z)\]
to be the full subcategory consisting of objects all of whose cohomology sheaves are
supported on the zero section $Z\subset \omega_Z$. Of course, when we say an object $E\in \Coh \omega_Z$  is supported on $Z$ we mean only that its reduced support is contained in $Z$; the scheme-theoretic support of $E$ will in general be some non-reduced fattening of $Z$, and $E$ will not be of the form $s_*(F)$ for any  $F\in \Coh Z$.

\subsection{The rolled-up helix algebra}

Let $(E_0,\cdots ,E_{n-1})$ be a simple collection in $\D$ and let
$(E_i)_{i\in \Z}$ be the helix it generates. The graded algebra
\[\bigoplus_{k\geq 0}\prod_{j-i=k} \Hom_{\D}(E_i,E_j)\]
is a variant of what Bondal and Polishchuk called the helix algebra. It carries a natural $\Z$-action coming from the
isomorphisms
\[\tensor \omega_{Z}\colon \Hom_{\D}(E_i,E_j)\lra
\Hom_{\D}(E_{i-n},E_{j-n}).\] Define the rolled-up helix algebra to be the
invariant subalgebra
\[B=\bigg[\bigoplus_{k\geq 0}\prod_{j-i=k}\Hom_{\D}(E_i,E_j)\bigg]^{\Z}.\]
The degree zero part $B_0$ has a basis consisting of the
idempotents
\[e_i=\prod_{j\equiv i \ (n)}\id_{E_j}\in \prod_j\End_{\D}(E_j),\]
and there are corresponding simple right $B$-modules $T_i$ defined by
\[\dim_{\C} (T_j e_i)=\delta_{ij}.\]
In contrast to the situation with the finite-dimensional algebras considered
in the last section these will not be the only simple $B$-modules.

\begin{prop}
\label{equivalence}
Let $(E_0,\cdots, E_{n-1})$ be a simple collection on $\D$ and
let $B$ be the associated rolled-up helix algebra.
Then the functor 
\[ \Fw=\RHom_{\DD}^{\blob}\big(\bigoplus_{i=0}^{n-1} \pi^* E_i, -\big)\colon
\D^b(\Coh \omega_Z) \lra \D^b(\Mod B)\]
is an equivalence of categories.
 \end{prop}

\begin{pf}
Note that $\pi_* (\OO_{\omega_Z})=\bigoplus_{p\leq 0} \omega_Z^p$. The
adjunction $\pi^*\dashv \pi_*$ together with the projection formula shows
that for arbitrary objects $E$ and $F$ of $\D(Z)$
\[\Hom_{\DD}^k(\pi^*E,\pi^* F)=\bigoplus_{p\leq
0}\Hom_{\D}^k(E,F\tensor\omega_Z^p).\] Since $(E_0,\cdots ,E_{n-1})$ is a
simple collection, it follows that
\[\End^k_{\DD}\big(\bigoplus_{i=0}^{n-1} \pi^* E_i\big)=
\biggl\{\begin{array}{ll} B & \text{ if }k=0, \\ 0 &\text{
otherwise.}\end{array}\biggr.\]
One has to play around with the adjunction maps a little to see that the algebra structure is the one described above. Applying the adjunction $\pi^*\dashv \pi_*$
again shows that for any object $E\in \Dw$
\[\Hom_{\DD}^k(\pi^* E_i, E)=0\text{ for all }k\in \Z \implies \pi_*(E)=0.\]
But the functor $\pi_*$ is an exact functor on the category $\Coh(\omega_Z)$
and has no kernel, so this implies that $E\isom 0$. The statement then
follows from the general theory of derived Morita equivalence \cite{Ri}.
\end{pf}

Under the equivalence $\Fw$, the object
$\pi^* E_i$ is mapped to the projective module $P_i=e_i B$, and if $(F_0,
\cdots, F_{n-1})$ is the dual collection to $(E_0,\cdots,E_{n-1})$
as in Lemma \ref{baah}, then the
object \begin{equation}
\label{knick}
S_j=s_* (F_{n-1-j}[j])
\end{equation}
is mapped to the simple module $T_j$.

\begin{prop}
If $(E_0,\cdots, E_{n-1})$ is a simple collection in $\D$ then the corresponding rolled-up helix algebra $B$ is generated over $B_0$ by $B_1$ and is Koszul.
\end{prop}

\begin{pf}
This is entirely analogous to the proof of Proposition \ref{generated}. It is
basically a corollary of Bondal and Polishchuk's result Theorem \ref{BondalPolishchuk}.
\end{pf}

The graded algebra $B$ can naturally be viewed
as the path algebra of a quiver with relations.
The quiver has $n$ vertices $\{0,1,\cdots,n-1\}$ corresponding to the
idempotents $e_i\in B_0$. For each $1\leq i \leq n-1$ there are
\[d_i=\dim \Hom_{\D}(E_{i-1},E_{i})\] arrows from vertex $i-1$  to vertex $i$.
The only difference to the quivers considered in the last section
is that there are now \[d_0=\dim
\Hom_{\D}(E_{n-1},E_n)\] arrows from vertex $n-1$ to vertex $0$. Thus
the quiver is a cycle
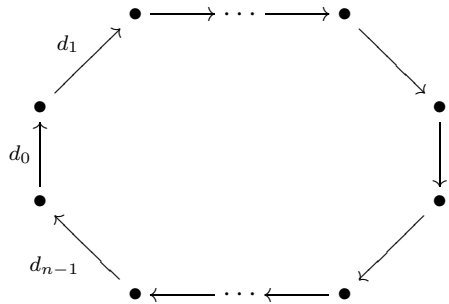
\[\xymatrix{ &  \bullet \ar[r] &\cdots \ar[r] & \bullet \ar[dr] \\
\bullet \ar[ur]^{d_1} &&&& \bullet \ar[d]
\\ \bullet \ar[u]^{d_{0}}
 &&&& \bullet \ar[dl] \\ & \bullet \ar[ul]^{d_{n-1}} &\cdots \ar[l] &\bullet \ar[l]}\]
As before, the Koszul property implies that the relations are quadratic.

\begin{example}
Set $Z=\PP^{n-1}$ and consider the diagonal action of the cyclic group $\Z_n$ on
affine space $\C^n$ with weights $\exp(2\pi i/n)$. The quotient
variety $X=\C^n/\Z_n$ has an isolated singularity; blowing it up gives the
variety $\omega_Z$; the resulting birational morphism contracts the
zero section $Z\subset\omega_Z$, and is a crepant resolution of
singularities.

The abelian category of $\Z_n$-equivariant coherent sheaves on $\C^n$ is
tautologically equivalent to the module category $\Mod R$ of the
corresponding skew group algebra  $R=\C[x_1,\cdots,x_n]*\Z_n$. We claim that
the ring $R$ is in fact isomorphic to the rolled-up helix algebra $B$ of the
helix $(\OO(i))_{i\in\Z}$ on $Z$, so that in this very special case, the
equivalence $\Fw$ can be thought of as an incarnation of the McKay
correspondence.

To prove the claim, note first that the degree zero part of both graded
algebras $B$ and $R$ is the same, namely a semisimple algebra spanned by
idempotents $e_0,\cdots, e_{n-1}$. Furthermore, for all $0\leq i\leq j\leq n-1$ there
are natural identifications
\[e_i B e_j= e_i R e_j =\C[x_1,\cdots, x_n]^{(j-i)},\] where the right hand
side is the space of polynomials of degree congruent to $j-i$ modulo $n$. It
is easy to check that the maps \[ e_i B e_j \tensor e_j B e_k\to e_i B
e_k,\qquad e_i R e_j \tensor e_j R e_k\to  e_i R e_k\] correspond to
multiplication of polynomials, and so the claim follows.
\end{example}

A right module $M$ over $B$ is said to be \emph{nilpotent} if there is some natural
number $n$ such that $M B_n=0$. Let $\NilMod B\subset \Mod B$ denote the
thick abelian subcategory consisting of nilpotent modules. Since any
module satisfying $M B_1=0$ is a direct sum of copies of the simple modules
$T_i$, one sees that $\NilMod B$ is a finite length category with simple
objects $T_0,\cdots, T_{n-1}$. In fact it is the smallest extension-closed
subcategory of $\Mod B$ containing each module $T_i$.

Let $\D_0^b(\Mod B)\subset\D^b(\Mod B)$
be the full subcategory consisting of objects whose cohomology modules are nilpotent.
It is not immediately clear whether this category can be identified with the derived category $\D^b(\NilMod B)$. A similar question arises as to whether $\Dw$ is the derived category of the subcategory of $\Coh \omega_Z$ consisting of sheaves supported on the zero section. But these questions will not be important for us.

\begin{lemma}
\label{small}
The equivalence
\[ \Fw\colon \D^b(\Coh \omega_Z) \lra \D^b(\Mod B)\]
of Proposition \ref{equivalence} restricts to give an equivalence
of full subcategories
\[\Fw\colon \Dw\lra\D_0^b(\Mod B).\]
\end{lemma}

\begin{pf}
This is immediate since $\Dw$ is the smallest full triangulated subcategory of
$\D$ containing the objects $S_j$ and $\D_0^b(\Mod B)$ is the smallest full
triangulated subcategory of $\D^b(\Mod B)$ containing the simple modules $T_j$.
\end{pf}

\subsection{Spherical collections}

In Section 3, rather than working directly with a given exceptional subcategory of $\D$, we worked with the corresponding set of projective objects, which formed an exceptional collection $(E_0,\cdots, E_{n-1})$. We then used the braid group action on exceptional collections to get a handle on the combinatorics of the exceptional subcategories. Of course, we could equally well have worked with the simple objects of a given exceptional subcategory, which are closely related to the dual exceptional collection $(F_0,\cdots, F_{n-1})$.

In the next subsection we shall be interested in certain finite length abelian subcategories of $\Dw$. Neither the projective nor the simple objects of these subcategories form exceptional collections. However, in this case, the simples are what Seidel and Thomas \cite{ST} called spherical objects, and together they form what we shall call a spherical collection. In this subsection we define an action of the group $B_n$ on the set of spherical collections in $\Dw$; this will be used in the next subsection to analyse the combinatorics of the corresponding subcategories of $\Dw$.

\smallskip

Let $n$ be the dimension of the variety $\omega_Z$. An object $S\in \Dw$ is
\emph{spherical} if
\[
\Hom^k_{\Dw}(S,S)=\biggl\{\begin{array}{ll} \C & \text{ if }k=0 \text{ or }n,
\\ 0 &\text{ otherwise.}\end{array}\biggr.
\]
Since $\omega_Z$ has trivial canonical bundle, and any object $S\in \Dw$ has compact support, Serre duality gives an isomorphism of functors
\[ \Hom_{\Dw}(S,-)\isom \Hom_{\Dw}(-,S[n])^*.\]
 The following result then follows from constructions given in \cite{ST}.

\begin{prop}[Seidel, Thomas]\label{ST} If $S\in \Dw$ is spherical then there is an
auto--equivalence $\Phi_S\in \Aut \Dw$ such that for any $F\in
\Dw$ there is a triangle
\[\Hom_{\Dw}(S,F)\tensor S\lra F\lra \Phi_S(F).\]
Furthermore, $\Phi_{S[1]}\isom\Phi_S$, and one has relations
\[\Phi_{S_1}\circ \Phi_{S_2} \circ \Phi_{S_1}^{-1}\isom\Phi_{\Phi_{S_1}(S_2)},\]
for any pair of spherical objects $S_1,S_2\in\Dw$. \qed
\end{prop}

The autoequivalences $\Phi_S$ associated to spherical objects are often called \emph{twist} functors. A ready supply of spherical objects on $\omega_Z$ is obtained by extending
exceptional objects on $Z\subset \omega_Z$ by zero.

\begin{lemma}
\label{msri}
If $E\in \D$ is exceptional then $s_* E\in \Dw$ is spherical. More generally,
if $E$ and $F$ are objects of $\D$ satisfying
$\Hom_{\D}^k(E,F)=0=\Hom_{\D}^k(F,E)$
for all $k\neq 0$, then one has
\[\Hom_{\Dw}^{\blob}(s_* E, s_* F)=\Hom_{\D}(E,F)\oplus \Hom_{\D}(F,E)^*[-n].\]
\end{lemma}

\begin{pf}
If $s\colon Z\into Y$ is the inclusion of a smooth projective subvariety $Z$ in a smooth quasi-projective variety $Y$ then a standard calculation shows that for any pair of objects $E$ and $F$ of $\D^b(\Coh Z)$ there is a spectral sequence
\[\Hom^p_Z(E, F\tensor \wedge^q \mathcal{N})\implies \Hom^{p+q}_Y(s_* E, s_* F),\]
where $\mathcal{N}$ is the normal bundle of $Z$ in $Y$. Our result follows by taking $Y$ to be the total space of $\omega_Z$, so that $\mathcal{N}=\omega_Z$, and 
computing $\Hom_Z^{\blob}(E,F\tensor\omega_Z)$ using Serre duality.
\end{pf}

Define a spherical collection of length $n$ in $\Dw$ to be an ordered
collection of spherical objects $(S_0,\cdots ,S_{n-1})$. The following action of the group $B_n$ should be compared with the action of $A_n$ on exceptional collections described in Theorem \ref{exceptionalaction}. The formula given here is justified by Proposition \ref{hannah} below.

\begin{lemma}
\label{sphericalaction} The group $B_n$ acts on the set of length $n$
spherical collections in $\Dw$. The generator $r$ acts by
\[r(S_0,S_1,\cdots ,S_{n-1})=(S_{n-1},S_0,\cdots ,S_{n-2}),\]
and for $1\leq i\leq n-1$, the generator $\tau_i$ acts by
\[\tau_i(S_0,\cdots ,S_{n-1})=(S_0,\cdots,
S_{i-2},S_i[-1],\Phi_{S_i}(S_{i-1}),S_{i+1},\cdots ,S_n).\]
\end{lemma}

\begin{pf}
Note first that it is not necessary to define the action of $\tau_0$ since
$\tau_0=r^{-1}\tau_1 r$. Assume $n\geq 3$ and consider the relation
$\tau_1\tau_2\tau_1=\tau_2\tau_1\tau_2$. This is easy to check directly using
the relations of Lemma \ref{ST}; up to isomorphism both sides take the
spherical collection $(S_0,\cdots,S_{n-1})$ to the collection
\[(S_2[-2], \Phi_{S_2} (S_1)[-1], \Phi_{S_2}\Phi_{S_1} (S_0), S_3, \cdots, S_{n-1}).\]
The other relations are either obvious or follow from this by conjugating by
$r$.
\end{pf}

Note that the group of exact autoequivalences of $\Dw$ acts on the set of spherical
collections in the obvious way: if $\Phi\in \Aut \Dw$ is an exact
autoequivalence, and $(S_0,\cdots,S_{n-1})$ is  a spherical collection, then
\[\Phi(S_0,\cdots,S_{n-1})=(\Phi(S_0),\cdots, \Phi(S_{n-1})).\]
The elements
$\alpha_i=r^i(\tau_1\cdots \tau_{n-1})r^{-(i+1)}\in B_n$ defined in Lemma \ref{exactsequence} act on spherical collections by autoequivalences.

\begin{lemma}
\label{done}
If $(S_0,\cdots, S_{n-1})$ is a spherical collection in $\Dw$ then
\[\alpha_i(S_0,\cdots, S_{n-1})=\Phi_{S_i}(S_0,\cdots, S_{n-1}) \text{ for }0\leq i\leq n-1.\]
\end{lemma}

\begin{pf}
This is a simple computation using the definition of the action of $B_n$ in Lemma \ref{sphericalaction}. We leave the details to the reader.
\end{pf}

\subsection{T-structures and tilting}

Let $(E_0,\cdots,E_{n-1})$ be a simple collection in $\D$ and let $B$ be the corresponding rolled-up helix algebra.
The standard t-structure on $\D^b(\Mod B)$ induces one on $\D_0(\Mod B)$ in the obvious way, and pulling this back using the equivalence
\[\Fw\colon \Dw\lra\D_0(\Mod B)\]
 of Lemma \ref{small} gives a bounded t-structure on $\Dw$ whose heart
\[\Aw(E_0,\cdots,E_{n-1})\subset\Dw\] is equivalent to $\NilMod B$. Let us call
the subcategories of $\Dw$ obtained from simple collections in $\D$ in  this way \emph{exceptional}.

We shall also define a \emph{quivery} subcategory of $\Dw$ to be one of the form $\Phi(\Aw)\subset \Dw$ for some autoequivalence
$\Phi\in\Aut \Dw$ and some exceptional subcategory $\Aw\subset \D$. Note that the analogous definition in the last section would have given nothing new, since if $\Phi\in \Aut \D$ and $\A\subset \D$ is an exceptional subcategory corresponding to the exceptional collection $(E_0,\cdots, E_{n-1})$ then $\Phi(\A)\subset \D$ is the exceptional subcategory corresponding to the exceptional collection $\Phi(E_0,\cdots, E_{n-1})$.  

Any quivery subcategory of $\Dw$ is a finite length abelian category
with $n$ simple objects $S_0,\cdots,S_{n-1}$. By (\ref{knick}) and Lemma \ref{msri} these simple objects are spherical. They have a canonical cyclic
ordering in which
\begin{equation}
\label{stones}
\Hom^k_{\Dw}(S_i,S_j)=0\text{ unless } 0\leq k\leq n\text{ and
}i-j\equiv k\mod n.
\end{equation}
If $\Aw=\Aw(E_0,\cdots,E_{n-1})$ is an exceptional
subcategory then its simples are given by (\ref{knick}), and
thus have a canonical ordering $(S_0,\cdots,S_{n-1})$ compatible with the
above cyclic ordering. One consequence of the following result is that this statement does not
extend in an obvious way to quivery subcategories.

\begin{prop}
\label{rotate} Let $(E_0,\cdots,E_{n-1})$ be a simple collection in $\D$, and let
$(E_i)_{i\in\Z}$ be the helix it generates. If $(S_0,\cdots,S_{n-1})$ are the
simples in the exceptional subcategory
$\Aw(E_0,\cdots,E_{n-1})$ with their canonical ordering, then
$\Phi_{S_{n-1}}(S_{n-1},S_0, \cdots, S_{n-2})$ are the simples in
$\Aw(E_{-1},E_0,\cdots,E_{n-2})$ with their canonical ordering.
\end{prop}

\begin{pf}
Let $(F_0,\cdots,F_{n-1})=\delta (E_0,\cdots,E_{n-1})$ be the dual
collection.
Since
\[(E_{-1},\cdots, E_{n-2})=(\sigma_1\cdots\sigma_{n-1})(E_0,\cdots,E_{n-1}),\]
Lemma \ref{conj} shows that the dual collection to $(E_{-1},\cdots,E_{n-2})$ is
\[(F_0',\cdots,F_{n-1}') =(\sigma_{n-1}\cdots\sigma_1)(F_0,\cdots,F_{n-1})
=(\L_{F_0}(F_1),\cdots,\L_{F_0}(F_{n-1}),F_0).
\]
Thus if $(S_0',\cdots, S_{n-1}')$ are the simples in $\Aw(E_{-1},\cdots,E_{n-2})$ with their canonical ordering, then
\[S'_0=s_* F_0\quad\text{and}\quad S_j'= s_* (\L_{F_0} F_{n-j} [j])\quad\text{ for } 1\leq j\leq n-1.\]
For each $1\leq j\leq n-1$, pushing forward the definition of a mutation and using Lemma \ref{msri} gives a triangle
\[ s_* (\L_{F_0} F_{n-j}[j-1]) \lra \Hom_{\Dw}(s_* F_0,s_* F_{n-j})\tensor s_* (F_0 [j-1])
\lra s_* (F_{n-j}[j-1]).\]
Rotating the triangle and using (\ref{knick})
we can reinterpret this as a triangle
\[ \Hom_{\Dw}(S_{n-1},S_{j-1})\tensor S_{n-1}\lra S_{j-1} \lra s_* (\L_{F_0} F_{n-j}[j])\]
From the definition of the twist functor $\Phi_{S_{n-1}}$ it follows that $S'_j=\Phi_{S_{n-1}}(S_{j-1})$ for $1\leq j\leq n-1$. Finally, any spherical object $S\in \Dw$ satisfies $\Phi_S(S)=S[1-n]$. Applying this to $S_{n-1}$ shows that $S'_0=\Phi_{S_{n-1}}(S_{n-1})$ which completes the proof.
\end{pf}

An ordered quivery subcategory of $\Dw$ is defined to be a quivery
subcategory together with an ordering of its simple objects compatible with
the canonical cyclic ordering. Note that an ordered quivery subcategory
determines and is determined by the corresponding spherical collection
$(S_0,\cdots,S_{n-1})$.

\begin{prop}
\label{hannah}
Suppose $(S_0,\cdots,S_{n-1})$ are the ordered
simples of an ordered quivery subcategory $\Aw\subset \Dw$. Then  for any
$0\leq i\leq n-1$ the tilted subcategory $\L_{S_i}
\Aw\subset\Dw$ is a quivery subcategory, and its simple objects
with their canonical cyclic order are given by the
spherical collection $\tau_i(S_0,\cdots,S_{n-1})$.
\end{prop}

\begin{pf}
By applying a power of $r$ to the spherical collection $(S_0,\cdots, S_{n-1})$ and thus changing the ordering of the simples we can assume that the simple we tilt at is $S_1$, or in other words, we can take $i=1$. Furthermore,
it is easy to see that we can apply an
autoequivalence of $\Dw$ without affecting the hypotheses or the conclusion of the Proposition. Thus, we
may assume that \[\Aw=\Aw(E_0,\cdots,E_{n-1})\] is an exceptional subcategory,
and using Proposition \ref{rotate}, we may assume further that $(S_0,\cdots,S_{n-1})$ have the corresponding canonical ordering.

Consider the mutated exceptional collection
\[(E_0',\cdots, E_{n-1}')=\sigma_1(E_0,\cdots, E_{n-1}).\]
We claim that the tilted subcategory
$\L_{S_1}(\Aw)$ is the exceptional subcategory
$\Aw'=\Aw(E_0',\cdots, E_{n-1}')$.
The proof of this goes in exactly the same way as that of Proposition \ref{helloo}.
The simple objects of $\Aw'$ with their canonical ordering are given by
\[(S_1[-1], S_1', S_2, \cdots, S_{n-1}),\]
where $S_1'$ is the universal extension
\[0\lra S_0\lra S_1'\lra \Ext^1_{\Aw}(S_1,S_0)\tensor S_1\lra 0.\]
As in Proposition \ref{helloo} it follows that $\Aw'=\L_{S_1}(\Aw)$.
But by the defintion of the twist functor $S_1'=\Phi_{S_1}(S_0)$ so the result follows.
\end{pf}

Combining this result with Lemma
\ref{sphericalaction} gives our main theorem.

\begin{thm}
\label{two} There is an action of the group $B_n$ on
the set of ordered quivery subcategories of $\Dw$. For each integer
$0\leq i\leq n-1$ the element $\tau_i$ acts on the underlying abelian subcategories by
tilting at the $i$th simple object. \qed
\end{thm}

We conclude this section with a remark concerning the exact sequence
\[1\lra F_n\lra B_n\lRa{h} A_n/\langle\gamma\rangle\lra 1\]
 of Lemma \ref{exactsequence}.
Consider an ordered quivery subcategory $\Aw_1\subset \Dw$ and its image $\Aw_2=\tau (\Aw_1)$  under the
action of an element $\tau\in B_n$. Using Proposition \ref{rotate} we can find exceptional subcategories $\Aw'_1$ and $\Aw'_2$ of $\Dw$ such that each subcategory $\Aw_i$, with the chosen ordering of its simples, is related to the corresponding exceptional subcategory $\Aw'_i$, with the canonical ordering of its simples, by an autoequivalence $\Phi\in \Aut \Dw$. Then
the two exceptional collections defining $\Aw'_1$ and $\Aw'_2$ are related by the action of some element
of the coset $h(\tau)$ in $A_n$. We shall not need this fact in what follows and we leave the proof to the reader.


\section{The case $Z=\PP^2$}

In this section we study in more detail the case when $Z=\PP^2$ is the
projective plane. 
Thus $\D$ denotes the derived category $\D^b(\Coh \PP^2)$ and $\Dw$ denotes the full subcategory of $\D^b(\Coh \omega_{\PP^2})$ consisting of objects whose cohomology sheaves are supported on the zero section. Note that in this case $\omega_Z$ is the line bundle $\OO(-3)$. An exceptional collection of length three will be called an exceptional triple.

\subsection{Markov triples}
Exceptional collections on $\PP^2$ were studied in detail by 
Gorodentsev and Rudakov \cite{Go, GR}. They discovered a connection between exceptional
triples and a certain Diophantine equation called the Markov
equation.

\begin{defn}
A Markov triple is an ordered triple of positive integers $(a,b,c)$
satisfying the equation
\[ a^2 + b^2 + c^2 = abc\]
The set of Markov triples will be denoted $\Mar$.
\end{defn}

A good proof of the following result is given by Bondal and Polishchuk \cite[Example 3.2]{BP}.

\begin{prop}(Gorodentsev, Rudakov)
\label{marv}
If $(E_0,E_1,E_2)$ is a strong exceptional triple in $\D$, then the positive
integers $(a,b,c)$ defined by
\[a=\dim \Hom_{\D}(E_0,E_1),\quad b=\dim \Hom_{\D}(E_1,E_2),\quad c=\dim
\Hom_{\D}(E_0,E_2)\] form a Markov triple.\qed
\end{prop}

It turns out that the space $\Mar$ carries a natural action of the group
$\PSL(2,\Z)$. Recall that
\[\PSL(2,\Z)=\Z_3 *\Z_2 =\big\langle w,v: w^3=v^2=1 \big\rangle\]
where $w$, $v$ and $u=wv$ can be represented by the matrices
\[u=\mat{1}{0}{1}{1},\qquad
v=\mat{0}{-1}{1}{0},\qquad w=\mat{0}{1}{-1}{1},\] respectively.
Define an action of $\PSL(2,\Z)$ on the set $\Mar$ of Markov triples by the
operations
\[ w \colon (a,b,c) \mapsto (c,a,b),\qquad v\colon (a,b,c)\mapsto
(b,a,ab-c).\] The following result is due to Markov. For the readers
convenience, and since we could not find the exact statement in the
literature, we include a proof, essentially lifted from Cassels \cite{Cas}.

\begin{prop}
\label{freemarkov}
The induced action of the normal subgroup \[\Gamma^3=\Z_2 * \Z_2 * \Z_2=\big\langle v,w^{-1}vw,wvw^{-1}\big\rangle\subset \PSL(2,\Z)\]
of index three
 on the set
$\Mar$ of Markov triples is free and transitive.
\end{prop}

\begin{pf}
For the description of $\Gamma^3$ as a free product see \cite[Theorem 1.3.2]{Ra}.
Define the weight of a Markov triple $(a,b,c)$ to be the product $abc$. It is
enough to show that for any Markov triple $(a,b,c)\neq (3,3,3)$, exactly one
of the triples
\begin{equation}
\label{markov} (b,a,ab-c), \qquad (c,ac-b,a),\qquad (bc-a,c,b),
\end{equation}
has smaller weight. Indeed, this implies that for each $(a,b,c)\in\Mar$ there
is a unique element of $\Gamma^3$ taking $(a,b,c)$ to $(3,3,3)$.

To prove the claim, first suppose that $a,b,c$ are not all distinct. Without
loss of generality assume that $b=c$. Then $a^2+2b^2=ab^2$ and $b$ divides
$a$. Writing $a=db$ it follows that $d$ divides $2$, and the only
possibilities are $(3,3,3)$ and $(6,3,3)$, for which the claim can be checked
directly.

Thus we can assume that $a,b,c$ are distinct, and without loss of generality
we can take $a>b>c$. Note that
\[c(ab-c)=a^2 + b^2.\]
Since $a^2+b^2> c^2$ it follows that $ab-c>c$ so that the first triple of
$(\ref{markov})$ has larger weight than $(a,b,c)$. The same argument applies
to the second triple.

Reducing modulo three shows that each of $a$, $b$ and $c$ is divisible by
three. Consider the quadratic function
\[f(t)=t^2+b^2+c^2-tbc.\]
This has roots $a$ and $bc-a$. Since $f(b)<3 b^2-b^2c\leq 0$ it follows that
$b$ lies between these two roots, and hence $bc-a<a$. Thus the third triple
of (\ref{markov}) has smaller weight than $(a,b,c)$.
\end{pf}

It is natural to view the points of $\Mar$ as the vertices of a graph, with
two triples being connected by an edge if they are obtained one from the
other by one of the generators $v$, $w^{-1}vw$, $wvw^{-1}$ of
$\Gamma^3$. Clearly, the resulting graph is a tree, and is just the Cayley graph of
$\Gamma^3$ with respect to the given generators. This tree is known as the Markov tree;
it is perhaps most natural to draw it in the hyperbolic plane because $\PSL(2,\RR)$ is the corresponding group
of isometries.

\subsection{T-structures on $\D$}

Gorodentsev and Rudakov showed that if $(E_0,E_2,E_2)$ is an exceptional triple in $\D$ then each object $E_i$ is a shift of a locally-free sheaf on $\PP^2$. They also proved the following transitivity result.

\begin{prop}(Gorodentsev, Rudakov)
\label{GRR}
The braid group $A_3$ acts transitively on the set of exceptional triples of sheaves on $\PP^2$.\qed
\end{prop}

It follows that an exceptional triples in $\D$ is simple if and only if it consists of sheaves. 
Let $\Str(\PP^2)$ denote the set of exceptional subcategories of $\D$.
We consider $\Str(\PP^2)$ as a graph in which two  subcategories are linked by an edge if they are related by a tilt at a simple. Proposition \ref{GRR} implies that the connected components of the graph $\Str(\PP^2)$ are indexed by the integers, and all components are isomorphic.

It is well known that there is a short exact sequence
\[ 1\lra \Z\lra A_3 \lRa{f} \PSL(2,\Z)\lra 1,\]
where the map $f$ takes the generators $\sigma_1,\sigma_2$ of $B_3$ to the
elements $w^{-1}v$ and $vw^{-1}$ of $\PSL(2,\Z)$ respectively. The kernel of $f$ is generated by the element $\gamma=(\sigma_1\sigma_2)^3$.

We can define a map \[T\colon \Str(\PP^2)\lra \Mar\]
by sending an exceptional subcategory $\A\subset \D$ with ordered simples
$(S_0,S_1,S_2)$ to the triple of positive integers
\[a=\dim \Hom^1_{\D}(S_1,S_0),\quad b=\dim \Hom^1_{\D}(S_2,S_1),
\quad c=\dim \Hom^2_{\D} (S_2,S_0).\]
These form a Markov triple by Proposition \ref{marv} since the Homs between the simples are just the Homs between the objects of the exceptional collection dual to the one defining $\A$.

\begin{thm}
The action of the group $A_3$ on the set $\Str(\PP^2)$ of exceptional subcategories of
$\D=\D(\PP^2)$  is free. The map $T$ is equivariant, which is to say
\[T(\sigma \A)=f(\sigma) T(\A),\]
for any exceptional subcategory $\A\subset \D$ and any element $\sigma\in A_3$.
Two subcategories lie in the same fibre of $T$ precisely if they are related by an autoequivalence of $\D$.
\end{thm}

\begin{pf}
First we show that $T$ is equivariant. Let $\A=\A(E_0,E_1,E_2)$
be an exceptional subcategory of $\D$. If \[(F_0,F_1,F_2)=\delta(E_0,E_1,E_2)\] is the dual collection,
then the simple objects of $\A$ with their canonical ordering are $(F_2,F_1[1],F_0[2])$.
If we apply $\sigma_1$ to $\A$ then by Lemma \ref{conj} the dual collection changes by $\sigma_2$. Thus the new simples are $(F_1,\L_{F_1}(F_2)[1],F_0[2])$.
Consider the defining triangle
\[\L_{F_1}(F_2)\lra \Hom_{\D}(F_1,F_2)\tensor F_1\lra F_2.\]
Applying the functor $\Hom_{\D}(-,F_1)$ immediately gives
\[\Hom_{\D}(\L_{F_1}(F_2),F_1)=\Hom_{\D}(F_1,F_2).\]
Applying the functor $\Hom_{\D}(F_0,-)$ and using the fact that the mutated collection is strong gives a short exact sequence
\[0\lra \Hom_{\D}(F_0,\L_{F_1}(F_2))\lra \Hom_{\D}(F_1,F_2)\tensor\Hom_{\D}(F_0,F_1)\lra \Hom_{\D}(F_0,F_2)\lra 0.\]
Thus if $T(\A)=(a,b,c)$ then \[T(\sigma_1(\A))=(a,ab-c,b)=(w^{-1}v) (a,b,c)=f(\sigma_1)T(\A).\]
A similar argument for $\sigma_2$ completes the proof of equivariance.

Next we show that the action of $A_3$ is free. Suppose an element
$\sigma\in A_3$ fixes an exceptional subcategory $\A\subset \D$.
Since
the action of $\PSL(2,\Z)$ on $\Mar$ is transitive we may assume that
$T(\A)=(3,3,3)$. By Proposition \ref{freemarkov}, the stabilizer subgroup of $(3,3,3)$ in $\PSL(2,\Z)$ is generated by $w$. Since $f(\zeta)=w$ and the kernel of $f$ is generated by $\zeta^3$ it follows that $\sigma=\zeta^k$ for some integer $k$.

By the relation (\ref{bondal}) the element $\gamma=\zeta^3$ acts on exceptional
collections by twisting by the anticanonical bundle. If $L$ is any ample line
bundle on $Z$ then the only objects of $\D$ satisfying $E\tensor L\isom E$ are
those supported in dimension zero, and these cannot be exceptional since they
are not rigid. Since the element $\sigma^k=\zeta^{3k}$ of $A_3$
fixes $\A$, and hence the exceptional objects which define it, it follows that $k=0$, which proves that the action is free.

For the last statement, note first that one implication is trivial since $T$ is defined in terms of dimensions of Hom spaces, and these are preserved by autoequivalences.
For the converse, observe that the action of $\Aut \D$ on $\Str (\PP^2)$
commutes with the action of $A_3$, so it will be enough to check that if two
exceptional subcategories $\A_1$ and $\A_2$ both lie over $(3,3,3)$ then they
differ by an autoequivalence. By Proposition \ref{GRR} the action of $A_3$ on
$\Str(\PP^2)$ is
transitive (up to shift) so we can assume that $\A_1=\A(\OO,\OO(1),\OO(2))$
and $\A_2=\sigma (\A)$ for some $\sigma\in A_3$. But
as above, $\sigma=\zeta^k$ for some integer $k$, and so
\[\A_2=\sigma(\A_1)=\A(\OO(k),\OO(k+1),\OO(k+2)),\]
which differs from $\A_1$ by tensoring with the line bundle $\OO(k)$.
\end{pf}

\subsection{T-structures on $\Dw$}
Consider now the corresponding picture for the category $\Dw$. The exact sequence of Proposition \ref{exactsequence}
takes the form
\[ 1\lra \Z*\Z*\Z \lra B_3 \lRa{g} \PSL(2,\Z)\lra 1,\] where the map $g$
is given by
\[g(r)=w, \quad g(\tau_i)=w^{i+1} v w^{1-i}\text { for }i\in \Z_3.\]
Let $\Ord(\PP^2)$ denote the set of ordered quivery abelian
subcategories of $\Dw$. We can define a map
\[T\colon \Ord(\PP^2)\lra \Mar\]
by sending a quivery subcategory with ordered simples
$(S_0,S_1,S_2)$ to the positive integers
\[a=\dim \Hom^1_{\Dw}(S_1,S_0), \quad b=\dim \Hom^1_{\Dw}(S_2,S_1), \quad
c=\dim \Hom^1_{\Dw}(S_0,S_2).\]
Once again, these integers form a Markov triple because by (\ref{knick}) and Lemma \ref{msri} the
Hom spaces coincide with Hom spaces between the objects of an exceptional collection.

\begin{thm}The action of the group $B_3$ on the set $\Ord(\PP^2)$ of ordered quivery subcategories of
$\Dw$  is free. The map $T$ is equivariant, which is to say
\[T(\tau \Aw)=g(\tau) T(\Aw),\]
for any ordered quivery  subcategory $\Aw\subset \Dw$ and any element $\tau\in B_3$.
Two ordered subcategories lie in the same fibre of $T$ precisely if they are related by an autoequivalence of $\Dw$.
\end{thm}

\begin{pf}
The proof of the equivariance of $T$ is almost the same as the one given in the last subsection and we omit it. However the proof that the action of $B_3$ is free is somewhat more complicated in this case.
Suppose an element $\tau\in B_3$ fixes an ordered quivery subcategory with
simples $(S_0,S_1,S_2)$. Since the action of $\PSL(2,\Z)$ on $\Mar$ is transitive, we can assume that $T(S_0,S_1,S_2)=(3,3,3)$. The stabilizer subgroup of $(3,3,3)$ in $\PSL(2,\Z)$ is generated by $w$, and $g(r)=w$, so for some integer $k$ the element $\tau r^k\in B_n$
lies in the kernel of the map $g$, which is freely generated by the elements $\alpha_0,\alpha_1,\alpha_2$ of Lemma \ref{exactsequence}.
Thus it will be enough to show that the subgroup $\Gamma\subset B_3$ generated by $\alpha_1$ and $r$ acts freely on the fibre
\[F=T^{-1}(3,3,3)\subset \Ord(\PP^2).\]

The Grothendieck group $\K(\Dw)$ is a rank three free abelian group. The
Euler form defines a skew-symmetric bilinear form on $\K(\Dw)$. Any autoequivalence of $\Dw$ induces an isometry of $\K(\Dw)$. The quotient of $\K(\Dw)$ by the
kernel of the Euler form is a rank two abelian group $\Lambda$ with an induced
non-degenerate skew-symmetric form. Any ordered quivery subcategory $\Aw\subset \Dw$ determines three ordered simples objects $(S_0,S_1,S_2)$ and hence a basis
$([S_0],[S_1],[S_2])$ of $\K(\Dw)$ and a basis $([S_0],[S_1])$ of $\Lambda$.

We claim that if $\Aw\subset \Dw$ is an ordered quivery subcategory of $\Dw$ lying in the fibre $F$, then so are $\alpha_1(\Aw)$ and $r(\Aw)$, and the corresponding bases of $\Lambda$ are related by the matrices
\[u^{3}=\mat{1}{0}{3}{1}\text{ and }w^{-1}=\mat{-1}{1}{-1}{0},\]
respectively. By Lemma \ref{done}, if the ordered simples of $\Aw$ are $(S_0,S_1,S_2)$ then the ordered simples 
of $\alpha_1(\Aw)$ are given by $\Phi_{S_1}(S_0,S_1,S_2)$. Since $\Aw$ lies in the fibre $F$ we have equalities in $\K(\Dw)$
\[[\Phi_{S_1}(S_0)]=[S_0]+3[S_1],\quad [\Phi_{S_1}(S_1)]=[S_1]\]
which gives the first matrix.
The fact that $\Aw$ lies in the fibre $F$ implies that the kernel of the Euler
form is generated by $[S_0]+[S_1]+[S_2]$. This means that
$[S_2]=-[S_0]-[S_1]$ in $\Lambda$ which gives the second matrix.

According to \cite[Theorems 1.7.4, 1.7.5 and Table 4]{Ra}, the elements $u^3, w u^3 w^{-1}$ and $w^{-1} u^3 w$ freely generate the
normal subgroup
\[\Gamma(3)=\Z * \Z * \Z=\big\langle u^3, w^{-1}u^3 w, w u^3 w^{-1}\rangle
\subset\PSL(2,\Z),\]
and this group does not contain the
elements $w^{\pm 1}$, so it follows that $\Gamma$ acts freely on $F$.

Finally we have to prove that any two ordered quivery subcategories $\Aw_1$, $\Aw_2$ lying over $(3,3,3)$ differ by an autoequivalence. Using Lemma \ref{rotate} we can assume that the two subcategories are in fact exceptional and that the simples have the corresponding canonical ordering. Thus by Proposition \ref{GRR}, we can take $\Aw_1=\Aw(\OO,\OO(1),\OO(2))$ and $\Aw_2=\tau \Aw_1$ for some $\tau\in B_3$. As above, it follows that for some integer $i$ the element $\tau r^i$ lies in the kernel of $g$. But the kernel of $g$ acts by autoequivalences, and by Proposition \ref{rotate},
applying $r^i\Aw_1$ differs from $\Aw_1$ by an autoequivalence, 
so the result follows.  
\end{pf}

\bigskip

\noindent Department of Pure Mathematics,
University of Sheffield,
Hicks Building, Hounslow Road, Sheffield, S3 7RH, UK.

\smallskip

\noindent email: {\tt t.bridgeland@sheffield.ac.uk}

\end{document}